\documentstyle[12pt]{article}
\let\Bbb\relax
\newfont{\Bb }{msbm10 scaled 1000} 
\newfont{\Bbb}{msbm10 scaled 1200}  
\newfam\eufam%
\font\euzw=eufm10 scaled 1200%
\font\euac=eufm9%
\textfont\eufam=\euzw\scriptfont\eufam=\euac%
\def\fr{\fam\eufam\euzw}%
\newcommand{\Z}{{\Bbb Z}}
\newcommand{\R}{{\Bbb R}}

\newtheorem{thm}{Theorem}
\newtheorem{prop}{Proposition} 
\newtheorem{cor}{Corollary} 
\newtheorem{lemma}{Lemma} 
\newtheorem{dof}{Definition} 
\begin{document}
{\small 
\rightline{dg-ga/9704002}
\rightline{MPI 97-29}
}
\vskip 0.5 true cm
\begin{center} 
{\LARGE Killing spinors are Killing vector fields in \smallskip
Riemannian Supergeometry}
\vskip 0.7 true cm 
{\large D.V.\ Alekseevsky $^a$\footnote{Supported by 
Max-Planck-Institut f\"ur Mathematik (Bonn). e-mail: aba@math.genebee.msu.su},
V.\ Cort\'es $^c$\footnote{Supported by AvH Foundation, MSRI 
(Berkeley)
and  SFB 256 (Bonn University). e-mail: 
vicente@math.uni-bonn.de},
C.\ Devchand $^d$\footnote{e-mail: devchand@ictp.trieste.it},
U.\ Semmelmann $^s$\footnote{Supported by Graduiertenkolleg 
Bonn. e-mail: uwe@rhein.iam.uni-bonn.de}}\end{center}
\noindent
{\small
$^a$ Sophus Lie Center,
Gen.\ Antonova 2, kv.\ 99,
117279 Moscow,
Russian Federation\\
$^{c,s}$ Mathematisches Institut der Universit\"at Bonn,
Beringstr. 1,
53115 Bonn, Germany\\
$^d$ Max-Planck-Institut f\"ur Gravitationsphysik, Schlaatzweg 1, 
14473 Potsdam, Germany
}
\vskip 0.5 true cm
\begin{abstract} 
A supermanifold $M$ is canonically associated to any pseudo Riemannian 
spin manifold $(M_0,g_0)$. Extending the metric $g_0$ to a field $g$ 
of bilinear forms  $g(p)$ on $T_pM$, $p\in M_0$,  the pseudo Riemannian
supergeometry of  $(M,g)$ is formulated as $G$-structure on $M$, 
where $G$ is a supergroup with even part $G_0\cong Spin(k,l)$; 
$(k,l)$ the signature of $(M_0,g_0)$. Killing vector fields on $(M,g)$  
are, by definition, infinitesimal automorphisms of this $G$-structure.
For every spinor field $s$ there exists a corresponding odd vector field 
$X_s$ on $M$. Our main result is that $X_s$ is a Killing vector field 
on $(M,g)$ if and only if $s$ is a twistor spinor. In particular, 
any Killing spinor $s$ defines a Killing vector field $X_s$. 
\end{abstract} 

\section{Introduction to supergeometry} 
First we introduce the supergeometric language which is needed
to formulate the main result of the paper. Standard references 
on supergeometry are \cite{M}, \cite{L} and \cite{K}.

\bigskip\noindent
1.1 {\bf Supermanifold.} We consider pairs $(M_0,{\cal A})$, 
where $M_0$ is a $C^{\infty}$-manifold and 
${\cal A} = {\cal A}_0 +  {\cal A}_1$ is a sheaf of 
$\mbox{\Z}_2$-graded {\R}-algebras; $\dim M_0 = m$.

\medskip\noindent
{\bf Example 1:} We denote by ${\cal C}^{\infty}_{M_0}$ the
sheaf of (smooth) functions of $M_0$. It associates to an open set
$U\subset M_0$ the algebra ${\cal C}^{\infty}_{M_0}(U) = C^{\infty}(U)$
of smooth functions on $U$. Let $E$ be a (smooth) 
vector bundle over $M_0$
and ${\cal E}$ the corresponding locally free sheaf of 
${\cal C}^{\infty}_{M_0}$-modules: ${\cal E}$  associates to an open set
$U\subset M_0$ the $C^{\infty}(U)$-module  
${\cal E} (U) = \Gamma (U,E)$ of sections of $E$ over $U$. 
Conversely, any locally free sheaf 
${\cal E}$ of 
${\cal C}^{\infty}_{M_0}$-modules defines a vector bundle 
$E\rightarrow M_0$. The exterior sheaf 
$\wedge {\cal E} = 
\wedge^{ev} {\cal E} + \wedge^{odd} {\cal E}$ is a sheaf of 
$\mbox{\Z}_2$-graded {\R}-algebras on $M_0$.  

\begin{dof} \label{supermfDef} The pair 
$M = (M_0, {\cal A})$ is called a (differentiable)
supermanifold of dimension $m|n$ over $M_0$ if for all $p\in M_0$
there exists an open neighborhood $U\ni p$ and a rank 
$n$ free sheaf
${\cal E}_U$ of ${\cal C}^{\infty}_{U}$-modules over 
$U$ such that 
${{\cal A}|}_U \cong \wedge {\cal E}_U$ (as sheaves of 
$\mbox{\Z}_2$-graded {\R}-algebras). The (local) sections of 
$\cal A$ are called (local) functions
on $M$. 
\end{dof}

\noindent
{}From  Def.\ \ref{supermfDef} it follows that there exists a canonical 
epimorphism $\epsilon : {\cal A} \rightarrow {\cal C}^{\infty}_{M_0}$, 
which is called the {\bf evaluation map}. Its kernel is the ideal
$\cal J$ generated by ${\cal A}_1$: $\ker \epsilon = {\cal J} = 
\langle {\cal A}_1 \rangle = {\cal A}_1 + {\cal A}_1^2$. By the construction 
of Example 1 to any vector bundle $E\rightarrow M_0$ we have associated
a supermanifold $M(E) = (M_0, {\cal A} = \wedge {\cal E})$. In this case
the exact sequence 
\[ 0 \rightarrow {\cal J} =  \langle {\cal E}\rangle \rightarrow 
{\cal A} = \wedge {\cal E} \stackrel{\epsilon}{\rightarrow} 
{\cal C}^{\infty}_{M_0} \rightarrow 0 \]
of sheafs of $\mbox{\Z}_2$-graded 
{\R}-algebras has a canonical splitting ${\cal C}^{\infty}_{M_0} 
\hookrightarrow \wedge {\cal E} = {\cal C}^{\infty}_{M_0} + 
\langle {\cal E}\rangle$. 

Let $(x^1,\ldots , x^m)$ be local coordinates for $M_0$ defined on an open 
set $U\subset M_0$ such that ${{\cal A}|}_U \cong \wedge {\cal E}_U$, where
${\cal E}_U$ is a rank $n$ free sheaf of ${\cal C}^{\infty}_{U}$-modules,
cf.\ Def.~\ref{supermfDef}. Let $\theta_1, \ldots ,\theta_n$ be sections of 
${\cal E}_U$ trivializing the vector bundle $E_U$ associated to the sheaf
${\cal E}_U$. Note that $x^1,\ldots , x^m, \theta_1 ,\ldots , \theta_n$
can be considered as local functions on the supermanifold $M$. Moreover, 
any local function $f\in {\cal A}(U)$ is of the form 
\begin{equation} \label{thetaexpansionEqu} 
f = \sum_{\alpha \in \mbox{\Z}_2^n} f_{\alpha} 
(x^1,\ldots , x^m)\theta^{\alpha}\, , \quad 
f_{\alpha} (x^1,\ldots , x^m) \in C^{\infty}(U) 
= {\cal C}^{\infty}_{M_0} (U)\, ,
\end{equation}   
where $\theta^{\alpha} := \theta_1^{\alpha_1} \wedge \ldots \wedge 
\theta_n^{\alpha_n}$, $\alpha = (\alpha_1,\ldots ,\alpha_n)$. 

\begin{dof} The tupel $(x^i,\theta_j) = (x^1,\ldots , x^m, \theta_1,\ldots ,
\theta_n)$ is called a {\bf local coordinate system} for $M$ over $U$.
\end{dof} 

The evaluation map applied to a (local) function
$f = f(x^1, \ldots ,x^m, \theta_1, \ldots , \theta_n)$ with expansion
(\ref{thetaexpansionEqu}) is given by:
\[ \epsilon (f) = f(x^1, \ldots ,x^m,0,\ldots , 0) 
= f_{(0,\ldots , 0)}(x^1, \ldots ,x^m)\, .\]     

Let $M = (M_0,{\cal A})$ and $N = (N_0, {\cal B})$ be supermanifolds.

\begin{dof} \label{morphismDef} 
A {\bf morphism}  $\Phi : M\rightarrow N$ is a pair $\Phi = 
(\varphi,\phi)$, where $\varphi : M_0 \rightarrow N_0$ is a
differentiable map and 
$\phi : {\cal B} \rightarrow \varphi_{\ast} {\cal A}$ is a 
morphism of sheaves of $\mbox{\Z}_2$-graded {\R}-algebras. 
$\Phi$ is called an {\bf isomorphism} if $\varphi$ is a diffeomorphism
and $\phi$ is an isomorphism. An isomorphism $\Phi : M\rightarrow M$ is 
called {\bf automorphism} of $M$.   
\end{dof}

In local coordinate systems $(x^i,\theta_j)$ for $M$ and 
$(\tilde{x}^k,\tilde{\theta_l})$ for $N$ a morphism $\Phi$ is expressed by
$p$ even functions 
$\tilde{x}^k(x^1, \ldots ,x^m, \theta_1, \ldots , \theta_n)$, 
$k = 1, \ldots , p$, and odd $q$ functions 
$\tilde{\theta}_l(x^1, \ldots ,x^m, \theta_1, \ldots , \theta_n)$, 
$l = 1,\ldots , q$; where $(p,q) = \dim N$.
 
\bigskip\noindent
1.2 {\bf Tangent vector/vector field.} Let $M = (M_0,{\cal A})$ be a 
supermanifold. For any point $p\in M_0$ the evaluation map
$\epsilon : {\cal A} \rightarrow {\cal C}^{\infty}_{M_0}$ induces
an epimorphism ${\epsilon}_p: {\cal A}_p\rightarrow \mbox{\R}$,
${\epsilon}_p(f) := {\epsilon}(f)(p)$, where ${\cal A}_p$ denotes the stalk
of $\cal A$ at $p$. For $\alpha \in \mbox{\Z}_2 =\{ 0,1\}$ we define 
\[ (T_pM)_{\alpha} := \{ v: {\cal A}_p\rightarrow \mbox{\R} 
\quad \mbox{{\R}-linear}|v(fg) = v(f){\epsilon}_p(g) + 
{(-1)}^{\alpha \tilde{f}}
{\epsilon}_p(f) v(g) \}\, ,\]
where the equation is required for all $f,g\in {\cal A}_p$ of pure degree and 
$\tilde{f} \in \{ 0,1\}$ denotes the degree of $f$. 

\begin{dof} The {\bf tangent space} of $M$ at $p\in M_0$ is the 
$\mbox{\Z}_2$-graded vector space 
$T_pM = (T_pM)_0 + (T_pM)_1$.  The elements of $T_pM$ are called 
{\bf tangent vectors}. Any morphism $\Phi = (\varphi , \phi ) : 
M = (M_0, {\cal A}) \rightarrow N = (N_0, {\cal B})$ induces 
linear maps $d\Phi (p) : T_pM \rightarrow T_{\varphi (p)} N$, defined 
by $(d \Phi (p)  v)(f) := v(\phi_p(f))$, $p\in M_0$, $v\in T_pM$, 
$f\in {\cal B}_{\varphi (p)}$, where $\phi_p : {\cal B}_{\varphi (p)}
\rightarrow {\cal A}_p$ is the morphism of stalks associated to 
$\phi : {\cal B} \rightarrow {\varphi}_{\ast}{\cal A}$. The map
$d \Phi (p)$ is called the {\bf differential at} $p$ of $\Phi$. 
\end{dof}

\noindent
The sheaf ${Der} {\cal A}$ of derivations of $\cal A$ over {\R} is a sheaf
of $\mbox{\Z}_2$-graded $\cal A$-modules: ${Der} {\cal A} = ({Der} {\cal A})_0
+ ({Der} {\cal A})_1$, where 
\[ ({Der} {\cal A})_{\alpha} = \{ X: {\cal A} \rightarrow 
{\cal A} \quad   
\mbox{{\R}-linear}| X(fg) = X(f)g + {(-1)}^{\alpha \tilde{f}} fX(g) \}\, ,\] 
where the equation is required for all $f,g\in {\cal A}$ of pure degree.  

\begin{dof} The sheaf ${\cal T}_M = {Der} {\cal A}$ is called the {\bf tangent
sheaf} of $M = (M_0,{\cal A})$. The sections of ${\cal T}_M$ are called 
{\bf vector fields}. 
\end{dof}  

Any local coordinate system $(x^i,\theta_j)$ over $U$ gives rise 
to even vector fields $\frac{\partial}{\partial x^i}$ and odd 
vector fields $\frac{\partial}{\partial \theta_j}$ over $U$. 
The action of the vector fields 
$\frac{\partial}{\partial x^i}, \frac{\partial}{\partial \theta_j}$ on a 
function $f$ with expansion (\ref{thetaexpansionEqu}) is given 
by:
\[  \frac{\partial f}{\partial x^i}  =  \sum _{\alpha} 
\frac{\partial f_{\alpha}(x^1,\ldots , x^m)}{\partial x^i}\theta^{\alpha}
\, ,\] 
\[ \frac{\partial f}{\partial \theta_j} =  
\sum_{\alpha} \alpha_j (-1)^{\alpha_1 + \cdots + \alpha_{j-1}} 
f_{\alpha} (x^1,\ldots , x^m)\theta_1^{\alpha_1} \wedge \ldots \wedge 
\theta_j^{\alpha_j -1} 
\wedge \ldots \wedge \theta_n^{\alpha_n}\, .\]

Any vector field $X$ on $M$ over $U$ can be written as
\[ X = \sum_{i=1}^m X^i(x^1, \ldots ,x^m, \theta_1, \ldots , \theta_n)
\frac{\partial}{\partial x^i} + \sum_{j=1}^n 
Y^j(x^1, \ldots ,x^m, \theta_1, \ldots , \theta_n) 
\frac{\partial}{\partial \theta_j} \, ,\]
where $X^i, Y^j \in {\cal A}(U)$. 

If $\Phi = (\varphi , \phi ) : 
M = (M_0, {\cal A}) \rightarrow N = (N_0, {\cal B})$ is an isomorphism then 
$\varphi^{-1}$ and $\phi^{-1} : {\varphi}_{\ast} {\cal A} \rightarrow
{\cal B}$ exist and give rise to an isomorphism ${\cal A} \rightarrow
{\varphi}^{-1}_{\ast} {\cal B}$. The induced isomorphism between the 
corresponding sheaves of derivations is denoted by 
\[ d\Phi : {\cal T}_M \rightarrow {\varphi}^{-1}_{\ast} {\cal T}_N\]
and is called the {\bf differential} of $\Phi$. For any open 
$U \subset M_0$ the differential $d\Phi$ is expressed by an 
${\cal A}(U)$-linear map $d\Phi_U : {\cal T}_M(U)\rightarrow  
{\cal T}_N (\varphi (U))$, where the action of ${\cal A}(U)$ on 
${\cal T}_N (\varphi (U))$ is defined using the isomorphism 
${\cal A}(U) \stackrel{\sim}{\rightarrow} {\cal B} (\varphi (U))$ induced
by $\phi^{-1}$. 

Let $X$ be a vector field defined on some open set $U\subset M_0$ and 
$p\in U$. Then we can define the {\bf value} $X(p) \in T_pM$ of 
$X$ at $p$: 
\[ X(p) (f) := \epsilon_p(X(f))\, ,\quad f\in {\cal A}_p\, .\] 
However, unless $\dim M = m|n = m|0$, a vector field is not determined by
its values. 

Finally, we relate the tangent spaces and tangent sheaves of $M$ and $M_0$. 
Any even tangent vector $v\in (T_pM)_0$ annihilates the ideal ${\cal J} 
= \ker {\epsilon}$ 
in the exact sequence 
\begin{equation}\label{exactsequEqu} 
0 \rightarrow {\cal J} \rightarrow {\cal A} 
\stackrel{\epsilon}{\rightarrow} {\cal C}^{\infty}_{M_0}
\rightarrow 0 \end{equation}   
and hence defines a tangent vector to $M_0$. More explicitly, 
we define a map $\epsilon : T_pM \rightarrow T_pM_0$ by 
the equation
\[ \epsilon (v) (\epsilon(f)) = v_0 (f)\, , \]
where $v = v_0 + v_1 \in (T_pM)_0 + (T_pM)_1$, $f \in {\cal A}_p$ and 
$f \mapsto \epsilon (f)$ is the  evaluation map of stalks 
$\epsilon :{\cal A}_p\rightarrow ({\cal C}^{\infty}_{M_0})_p$. 

\begin{prop} \label{exactsequProp} There is a canonical exact sequence of  
$\mbox{\Z}_2$-graded
vector spaces:
\[ 0 \rightarrow (T_pM)_1  \rightarrow T_pM \stackrel{\epsilon}{\rightarrow}
T_pM_0 \rightarrow 0\, .\] 
In particular, $\epsilon$ induces a canonical isomorphism $(T_pM)_0 
\stackrel{\sim}{\rightarrow} T_pM_0$.  
\end{prop}

\noindent 
Similarly, on the level of tangent sheaves we define $\epsilon : 
{\cal T}_M \rightarrow {\cal T}_{M_0}$ by the equation
\[ \epsilon (X) (\epsilon (f)) = \epsilon (X_0(f))\, ,\]
where $X = X_0 +X_1 \in ({\cal T}_M (U))_0 + ({\cal T}_M (U))_1$, 
$f\in {\cal A}(U)$ and $U\subset M_0$ open. 

\begin{prop} There is a canonical exact sequence of  sheaves of 
$\cal A$-modules 
\begin{equation} 0 \rightarrow \ker \epsilon \rightarrow  
{\cal T}_M \stackrel{\epsilon}{\rightarrow} {\cal T}_{M_0} \rightarrow 0\, ,
\end{equation} 
where $\ker \epsilon = ({\cal T}_M)_1 + {\cal J}{\cal T}_M$. In particular,
there is the following exact sequence of $\cal A$-modules:
\[ 0 \rightarrow ({\cal J}{\cal T}_M)_0 \rightarrow  
({\cal T}_M)_0 \rightarrow {\cal T}_{M_0} \rightarrow 0\, .\]
\end{prop}

\bigskip\noindent
1.3 {\bf Frame/frame field/local coordinates.} 
\begin{dof} Let $V = V_0 + V_1$ be a $\mbox{\Z}_2$-graded vector space of 
{\bf rank} $m|n$, i.e.\ $\dim V_0 = m$ and $\dim V_1 = n$.  A {\bf basis} 
of $V$ is a tupel $(b_1,\ldots , b_{m+n})$ such that $(b_1,\ldots , b_m)$
is a basis of $V_0$ and $(b_{m+1},\ldots , b_{m+n})$ is a basis of $V_1$. 
Let $M = (M_0, {\cal A})$ be a supermanifold and $p\in M_0$. A {\bf frame} at 
$p$ is a basis of $T_pM$. A tupel $(X_1,\ldots X_{m+n})$ of vector fields 
defined on an open subset $U\subset M_0$ is called a {\bf frame field} if
$(X_1(p),\ldots X_{m+n}(p))$ is a frame at all points $p\in U$. We denote
by ${\cal F} (U)$ the set of all frame fields over $U$. The sheaf of sets
$U\mapsto {\cal F} (U)$ is called the {\bf sheaf of frame fields}.   
\end{dof}

\noindent 
Any local coordinate system $(x^i,\theta_j)$ over $U$ gives rise 
to the frame field $(\frac{\partial}{\partial x^i},
\frac{\partial}{\partial \theta_j})$ over $U$. 

\bigskip\noindent
1.4 {\bf Supergroup.}  Let ${\bf A} = {\bf A}_0 + 
{\bf A}_1$ be an associative $\mbox{\Z}_2$-graded {\R}-algebra with unit.
We will always assume that ${\bf A}$ is {\bf supercommutative}, i.e.\ 
$ab = (-1)^{\tilde{a}\tilde{b}} ba$ for all $a,b \in {\bf A}_0 \cup 
{\bf A}_1$. Under this assumption any left-${\bf A}$-module 
carries a canonical right-${\bf A}$-module structure and 
vice versa; so we will
simply speak of ${\bf A}$-modules.  
For any supermanifold $M = (M_0, {\cal A})$ the algebra
of functions ${\cal A}(M_0)$ is supercommutative, associative and has a unit.
   
For any set $\Sigma$ and non-negative integers $r,s$ we
denote by $Mat (r,s,\Sigma)$ the set of $r\times s$-matrices 
with entries in $\Sigma$
and put $Mat (r,\Sigma) := Mat (r,r,\Sigma)$. Any partition 
$(r = m + n, s = k + l)$ defines a $\mbox{\Z}_2$-grading on the 
${\bf A}$-module $V = Mat (r,s,${\bf A}$)$: 
\begin{eqnarray*} V_0 &=& \{  \left( \begin{array}{ll} 
A & B \\
C & D
\end{array} \right) | A \in Mat(m,k,{\bf A}_0),\: 
D\in  Mat(n,l,,{\bf A}_0),\\
 & & B\in Mat(m,l,{\bf A}_1),\: 
C \in Mat(n,k,{\bf A}_1) \}\, , \\  
V_1  &=& \{  \left( \begin{array}{ll} 
A & B \\
C & D
\end{array} \right) | A \in Mat(m,k,{\bf A}_1),\: 
D\in  Mat(n,l,,{\bf A}_1),\\
 & & B\in Mat(m,l,{\bf A}_0),\: 
C \in Mat(n,k,{\bf A}_0) \} \, . 
\end{eqnarray*}
The $\mbox{\Z}_2$-graded ${\bf A}$-module $V = V_0 + V_1$ is denoted by 
$Mat (m|n, k|l, {\bf A})$. Matrix multiplication turns $Mat (m|n,{\bf A}) := 
Mat (m|n, m|n, {\bf A})$ into an associative $\mbox{\Z}_2$-graded 
algebra with unit.  

\begin{dof} A {\bf super Lie bracket} on a $\mbox{\Z}_2$-graded vector space  
$V = V_0 + V_1$ is a bilinear map $[ \cdot , \cdot ] : V \times V 
\rightarrow V$ such that for all $x,y,z \in V_0 \cup V_1$ we 
have:  
\begin{enumerate} 
\item[i)] $\widetilde{[x,y]} = \tilde{x} +  \tilde{y}$,
\item[ii)] $[x,y] = - (-1)^{\tilde{x}\tilde{y}} [y,x]$ and 
\item[iii)] $[x,[y,z]] = [[x,y],z] + (-1)^{\tilde{x}\tilde{y}} [y,[x,z]]$. 
\end{enumerate} 
The pair $(V,[\cdot , \cdot ])$  is called a {\bf super Lie algebra}. 
\end{dof} 

\noindent
The supercommutator
\[ [X,Y] = XY - (-1)^{\tilde{X}\tilde{Y}} YX\, , \quad
X,Y \in Mat(m|n,{\bf A})_0 \cup Mat(m|n,{\bf A})_1 \]
defines a super Lie bracket on the $\mbox{\Z}_2$-graded vector space
$Mat (m|n,{\bf A})$. The super Lie algebra $(Mat (m|n,{\bf A}), 
[ \cdot ,\cdot ])$  is denoted by ${\fr gl}_{m|n}({\bf A})$. 
We put 
\[ GL_{m|n}({\bf A}) := \{ g \in Mat (m|n,{\bf A})_0| g \quad
\mbox{is invertible} \} \, .\]
Similarly, if $V$ is a  $\mbox{\Z}_2$-graded ${\bf A}$-module 
$End_{\bf A} (V)$ carries a canonical super Lie algebra 
structure, 
which is denoted by ${\fr gl}_{\bf A}(V)$. By definition 
$GL_{\bf A}(V)$ is the group of invertible elements of $End_{\bf A} (V)$. 
Finally, we will use the convention ${\fr gl}_{m|n} := 
{\fr gl}_{m|n} (\mbox{\R})$, ${\fr gl}(V) := {\fr gl}_{\mbox{\R}} (V)$,
$GL(V) := GL_{\mbox{\R}} (V)$. 

\begin{dof} 
A {\bf supergroup} $G$ is a contravariant functor $M\mapsto G(M)$ from 
the category of supermanifolds into the category of groups. 
Let $H$, $G$ be supergroups. We say that $H$ is 
a {super subgroup} of $G$ and write $H\subset G$ if $H(M) \subset 
G(M)$ is a subgroup and $H(\Phi ) = G(\Phi )|H(N)$ for all supermanifolds
$M$, $N$ and morphisms $\Phi : M\rightarrow N$. 
\end{dof}

\noindent
{\bf Example 2:} The {\bf general linear supergroup} $GL_{m|n}$ is the 
supergroup $M \rightarrow GL_{m|n} (M)$ obtained as composition of the 
following two functors:
\begin{enumerate} 
\item[i)] the contravariant functor $M = (M_0,{\cal A}) \rightarrow 
{\cal A}(M_0)$ from the category of supermanifolds into that of asssociative,
supercommutative algebras with unit, 
\item[ii)] the covariant functor ${\bf A} \rightarrow GL_{m|n} ({\bf A})$
from the category of associative, supercommutative algebras with unit into 
that that of groups. 
\end{enumerate}

\begin{dof} A {\bf linear super Lie algebra} $\fr g$ is a super Lie subalgebra
${\fr g} \subset {\fr gl}_{m|n}$ (for some $m|n$). A {\bf linear supergroup}
is a super subgroup $G\subset GL_{m|n}$ (for some $m|n$).  
\end{dof} 
   
\noindent
{\bf Example 3:} Let ${\fr g} \subset {\fr gl}_{m|n}$ be a linear super Lie 
algebra. For any associative, supercommutative algebra with unit $\bf A$ we 
can consider the super Lie algebra ${\fr g} \otimes {\bf A}\subset 
{\fr gl}_{m|n}({\bf A})$. Its even part ${\fr g}({\bf A}) := 
({\fr g} \otimes {\bf A})_0$ is a Lie algebra. If ${\bf A} = {\cal A} (M_0)$
is the algebra of functions of a supermanifold $M = (M_0,{\cal A})$ then
it is easy to see that the exponential series 
\[ \sum_{n=0}^{\infty} \frac{1}{n!}X^n\, , \quad X\in {Mat} (m|n,{\bf A})\, ,\]
converges (locally uniformly) to an element $\exp X \in GL_{m|n}({\bf A})$. 
Now let $G({\bf A})$
be the subgroup of $GL_{m|n}({\bf A})$ generated by $\exp {\fr g}({\bf A})$.
then the functor $M = (M_0,{\cal A})\mapsto G(M) := G({\cal A}(M_0))$ is a 
linear supergroup, which we denote by $\exp {\fr g}$.  

\bigskip\noindent
1.5 {\bf G-structure.} Let $M = (M_0,{\cal A})$ be a super manifold of 
$\dim M = m|n$. For any open subset $U\subset M_0$ we consider the 
supermanifold ${M|}_U := (U,{{\cal A}|}_U)$. The ge\-ne\-ral linear supergroup
$GL_{m|n}$ induces a sheaf ${\cal GL}_M$ of groups over $M_0$:
${\cal GL}_M (U) := GL_{m|n}({M|}_U) = GL_{m|n}({\cal A}(U))$, $U\subset M_0$
open. The group ${\cal GL}_M (U)$ acts naturally (from the right) on the set
${\cal F} (U)$ of frame fields over $U$. This action turns $\cal F$ into a 
sheaf of ${\cal GL}_M$-sets. Now let $G \subset GL_{m|n}$ be a linear
supergroup and $\cal G$ the corresponding sheaf of groups, i.e.\ 
${\cal G} (U) = G({M|}_U)$ for all open $U\subset M_0$. Since $\cal G$ is a 
sheaf of subgroups ${\cal G} \subset{\cal GL}_M$ the sheaf $\cal F$ of 
frame fields of $M$ is, in particular, a sheaf of $\cal G$-sets. 

\begin{dof}  Let $M = (M_0,{\cal A})$, $\dim M = m|n$, be a supermanifold
and $G \subset GL_{m|n}$ a linear supergroup. A {\bf G-structure} on $M$ is 
a sheaf ${\cal F}_G$ of ${\cal G}$-subsets ${\cal F}_G \subset {\cal F}$ such 
that for all $p\in M_0$ there exists an open neighborhood $U\ni p$ for which
${\cal G} (U)$ acts simply transitively on ${\cal F}_G(U)$. 
\end{dof} 

\noindent
{\bf Example 4:} For any supermanifold $M$, $\dim M = m|n$, the sheaf of 
frame fields $\cal F$ is a $GL_{m|n}$-structure. 

\bigskip\noindent
1.6 {\bf Automorphism of G-structure.} We denote by ${Aut}  (M)$ the group of 
all automorphisms of the supermanifold $M$, see Def.\ \ref{morphismDef}.
The differential $d\Phi :{\cal T}_M \rightarrow 
\varphi_{\ast}^{-1} {\cal T}_M$ of any $\Phi = (\varphi , \phi ) \in
{Aut} (M)$  induces an isomorphism  
${\cal F} \rightarrow \varphi_{\ast}^{-1} {\cal F}$, 
again denoted by $d\Phi$. Now let ${\cal F}_G\subset {\cal F}$ be a 
$G$-structure on $M$, for some linear supergroup $G \subset GL_{m|n}$. 
For simplicity we can assume that $G = \exp {\fr g}$ as in Example 3. 

\begin{dof} \label{autoDef} $\Phi = (\varphi , \phi ) \in {Aut} (M)$ 
is called an 
{\bf automorphism} of the $G$-structure ${\cal F}_G$ if 
$d\Phi {\cal F}_G \subset \varphi_{\ast}^{-1} {\cal F}_G$.  
\end{dof}
 
\noindent
We recall that any $p\in M_0$ has an open neighborhood $U$ such that 
${\cal G} (U)$ acts simply transitively on ${\cal F}_G (U)$. Such open
sets $U\subset M_0$ will be called {\bf small}. If $U\subset M_0$
is small then ${\cal F}_G (U) = E {\cal G} (U)$ for any frame field 
$E\in {\cal F}_G (U)$. Here the right-action of the group ${\cal G} (U)$
on ${\cal F}_G (U)$ is simply denoted by juxtaposition. 

\begin{prop} \label{autoProp} $\Phi\in {Aut} (M)$ is an automorphism of the 
$G$-structure 
${\cal F}_G$ iff 
\[ {d\Phi}_{U'} {E|}_{U'} \in {E|}_{\varphi (U')} {\cal G} (\varphi (U'))\]
for all small $U \subset M_0$, $E\in {\cal F}_G (U)$ and open $U' 
\subset U$ such that $\varphi (U') \subset U$. 
\end{prop} 

\noindent
For any open set $U \subset M_0$ the vector space ${\cal T}_M(U)^{m+n}$
of $(m+n)$-tupels of vector fields is naturally a right-module of the 
associative, $\mbox{\Z}_2$-graded algebra $Mat (m|n, {\cal A}(U))$. 
In particular, it is a right-module of the super Lie algebra 
${\fr g} \otimes {\cal A} (U) \subset  {\fr g}_{m|n} ({\cal A}(U))$. 
On the other hand, ${\cal T}_M(U)$ (and hence ${\cal T}_M(U)^{m+n}$) is 
naturally  a left-module for the super Lie algebra  ${\cal T}_M(U)$ of 
local vector fields. The action on ${\cal T}_M(U)$ is given by the adjoint 
representation, i.e.\ by the supercommutator ${ad}_XY = X \circ Y -
(-1)^{\tilde{X}\tilde{Y}} Y \circ X$, $X, Y\in {\cal T}_M(U)$ of pure 
degree. The corresponding action on ${\cal T}_M(U)^{m+n}$ is denoted
by $L_X$ (``Lie derivative''):
\[ L_XE := ([X,X_1], \ldots , [X, X_{m+n}])\, ,\quad
E = (X_1, \ldots ,X_{m+n}) \in {\cal T}_M(U)^{m+n}\, .\] 
Proposition \ref{autoProp} motivates the following definition.

\begin{dof} \label{infautoDef} A vector field $X$ on $M$ is an 
{\bf infinitesimal
automorphism} of the $G$-structure ${\cal F}_G$ if 
\[ L_{{X|}_U} {E|}_U \in {E|}_U ({\fr g} \otimes {\cal A} (U))\]  
for all small $U \subset M_0$, $E\in {\cal F}_G (U)$. 
\end{dof}

\section{Supergeometry associated to the spinor bundle}
2.1 {\bf The supermanifold M(S).} 
Let $(M_0,g_0)$ be a (smooth) pseudo Riemannian spinmanifold with spinor bundle
$S\rightarrow M_0$. The corresponding locally free sheaf of 
${\cal C}^{\infty}_{M_0}$-modules will be denoted by $\cal S$; 
${\cal S} (U) = \Gamma (U,S)$, $U\subset M_0$ open. To the vector bundle
$S\rightarrow M_0$ we associate the supermanifold $M: M(S) = (M_0, {\cal A} =
\wedge {\cal S})$. 

Consider the $\mbox{\Z}_2$-graded vector bundle $TM_0 + S^{\ast} 
\rightarrow M_0$ with even part $TM_0$ and odd part $S^{\ast}$.
\begin{prop} \label{T_pMProp} For any $p\in M_0$ there is a canonical 
isomorphism of 
$\mbox{\Z}_2$-graded vector 
spaces ${\iota}_p : T_pM_0 + S^{\ast}_p \stackrel{\sim}{\rightarrow}
T_pM$. 
\end{prop} 

\noindent
{\bf Proof:} We define ${\iota}_p^{-1} | (T_pM)_0 := \epsilon |(T_pM)_0$, see  
Prop.\ \ref{exactsequProp}.  Now it is sufficient to construct a
canonical isomorphism $S^{\ast}\stackrel{\sim}{\rightarrow} (T_pM)_1$. 
For any section $s\in \Gamma (U,S^{\ast})$ interior multiplication ${\iota}(s)$
by $s$ defines an odd derivation of the $\mbox{\Z}_2$-graded algebra 
${\cal A}(U) = \Gamma (U,\wedge S)$, i.e.\ a vector field $X_s := \iota (s)
\in {\cal T}_M(U)_1$. The value $X_s(p) \in (T_pM)_1$  depends only on 
$s(p) \in S^{\ast}_p$ and we can define $\iota_p (s(p)) 
:= X_s(p)$. $\Box$ 
  
Using the embedding ${\cal C}^{\infty}_{M_0} \hookrightarrow \wedge {\cal S}$,
we can consider ${\cal T}_M$ as a sheaf of ${\cal C}^{\infty}_{M_0}$-modules. 
Interior multiplication $s\mapsto \iota (s) = X_s$ defines a monomorphism 
$S^{\ast} \hookrightarrow ({\cal T}_M)_1$ of sheaves of 
${\cal C}^{\infty}_{M_0}$-modules. We want to extend this map to 
$\iota : {\cal T}_{M_0} + {\cal S}^{\ast} \rightarrow {\cal T}_M$. 
For a local vector field $X\in {\cal T}_{M_0} (U)$ on $M_0$ we put 
\[ {\iota} (X) \: := \: \nabla_X \in {\cal T}_M (U)_0 \, ,\]
where $\nabla$ is the canonical connection on $\wedge S$, i.e.\ 
the one induced by the Levi-Civita-connection on $(M_0,g_0)$. 

\begin{prop} \label{iotaProp} The map $\iota : {\cal T}_{M_0} + 
{\cal S}^{\ast} 
\rightarrow {\cal T}_M$  is a monomorphism of sheaves of $\mbox{\Z}_2$-graded
${\cal C}^{\infty}_{M_0}$-modules. Moreover, $\iota |{\cal T}_{M_0}$ 
defines a splitting of the sequence (\ref{exactsequEqu}), i.e.\ 
$\epsilon \circ \iota|{\cal T}_{M_0} = id$. 
\end{prop} 

\noindent 
Note that given any vector bundle $E$ and connection $D$ on $E$
we can canonically define $\iota_{E,D} : {\cal T}_{M_0} 
+ {\cal E}^{\ast} \hookrightarrow {\cal T}_M$, where $M = M(E)$ and 
$\cal E$ is the sheaf of local sections of $E$. In Prop.\ \ref{iotaProp}
we have $\iota = \iota_{S,\nabla}$.    

\bigskip\noindent
2.2 {\bf The coadjoint representation of the Poincar\'e super Lie algebras.}
Let $(V_0,\langle \cdot ,\cdot \rangle )$ be a pseudo Euclidean vector space
of signature $(k,l)$, $k+l = m$, and $V_1$ the spinor module of the group 
${Spin} (V_0)$, $n:= \dim V_1 = 2^{[\frac{m}{2}]}$. Put $V := V_0 + V_1$. 
The vector space ${\fr p}(V) := {\fr spin} (V_0) + V$ carries the structure 
of ${\fr spin} (V_0)$-module. We want to extend this structure to a 
super Lie bracket $[ \cdot , \cdot ]$ on ${\fr p}(V)$ which satisfies 
$[V_0 , V] = 0$ and  $[V_1 ,V_1] \subset V_0$. 
Such an extension is precisely given by a 
${Spin} (V_0)$-equivariant map $\pi : \vee^2 V_1 \rightarrow V_0$;
here $\vee^2$ denotes the symmetric square. 

\begin{dof} The structure of super Lie algebra defined on ${\fr p}(V)$ by the
map $\pi$ is called  a {\bf Poincar\'e super Lie algebra}. 
\end{dof} 

\noindent
We denote by $\rho : V_0  \rightarrow {End} (V_1)$ the (standard) Clifford 
multiplication. 

\begin{dof}\label{suitableDef} 
A bilinear form $\beta$ on the spinor module is called 
{\bf admissible} if 
\begin{enumerate} 
\item[1)] $\beta$ is symmetric or skew symmetric. We define the symmetry
$\sigma$ of $\beta$ to be $\sigma (\beta ) = +1$ in the first case and 
$\sigma (\beta ) = -1$ in the second. 
\item[2)] Clifford multiplication $\rho(v)$, $v\in V_0$, is either 
symmetric or skew symmetric. Accordingly, we define the type $\tau$ of
$\beta$ to be $\tau (\beta ) = \pm 1$. 
\end{enumerate} 
An admissible form $\beta$ is called {\bf suitable} if $\sigma (\beta ) 
\tau (\beta ) = + 1$. 
\end{dof}
Given a suitable bilinear form $\beta$ on $V_1$ we define  
$\pi = {\pi}_{\rho , \beta} : \vee^2 V_1 \rightarrow V_0$ by 
\begin{equation} \label{piEqu}  \langle \pi (s_1 \vee s_2 ), v\rangle \: 
= \: \beta (\rho (v) s_1 ,s_2 ) ,\quad s_1,s_2 \in V_1, \quad v\in V_0
\, .\end{equation}
The map $\pi$ is ${Spin} (V_0)$-equivariant. Hence it defines 
on the vector space ${\fr p}(V)$ the structure of Poincar\'e super Lie algebra.
The following theorem was proved in \cite{A-C}. 

\begin{thm} Any ${Spin} (V_0)$-equivariant map 
$\vee^2 V_1 \rightarrow V_0$ is a linear combination of maps 
$\pi_{\rho , \beta_i}$, $\beta_i$ suitable. 
\end{thm}
 
\noindent
All admissible bilinear forms on the spinor module were 
explicitly determined in \cite{A-C}. The spinor module 
carries a {\bf non-degenerate} suitable bilinear form
$\beta$ for all values of $m = k+l$ and $s = k-l$ except for 
$(m,s) = (5,7)$,  $(6,0)$, $(6,6)$ and $(7,7)$ $\pmod{(8,8)}$. 
Now we assume that a non-degenerate suitable bilinear form 
$\beta$ on  $V_1$ is given. 
The map $\pi = \pi_{\rho ,\beta}$ defines on ${\fr p}(V)$ the   
structure of Poincar\'e super Lie algebra such that 
$[V_1,V_1] = V_0$. 

Given a super Lie algebra $\fr g$
the {\bf coadjoint representation} ${ad}^{\ast} : {\fr g} \rightarrow 
{\fr gl} ({\fr g}^{\ast})$, $x\mapsto {ad}^{\ast}_x$, is defined by the
equation 
\[ {ad}^{\ast}_x (y^{\ast}) = - (-1)^{\tilde{x}\tilde{y^{\ast}}}y^{\ast} 
\circ {ad}_x\, ,\]
for $x\in {\fr g}$ and $y^{\ast} \in {\fr g}^{\ast}$ of pure degree. 

\begin{prop} \label{coadProp} The coadjoint representation of ${\fr p}(V)$ 
preserves the 
subspace $V^{\perp} = \{ x^{\ast} \in {\fr p}(V)^{\ast}| x^{\ast}(V) = 0\}
\subset {\fr p}(V)^{\ast}$ and hence induces a representation
$\alpha : {\fr p}(V) \rightarrow {\fr gl} (V^{\ast})$ on $V^{\ast} \cong 
{\fr p}(V)^{\ast}/V^{\perp}$. It has kernel $\ker \alpha = V_0$ and therefore
induces a faithful representation of the super Lie algebra ${\fr p}(V)/V_0$ on
$V^{\ast}$. 
\end{prop} 

\noindent
Once we choose a basis $b = (b_1,\ldots , b_{m+n})$ of $V^{\ast}$, we can 
identify $\alpha ({\fr p}(V)) \subset {\fr gl} (V^{\ast})$ with a subalgebra
$\alpha ({\fr p}(V))^b \subset {\fr gl}_{m|n}$, where $A\mapsto A^b$
denotes the isomorphism ${\fr gl}(V^{\ast}) \rightarrow {\fr gl}_{m|n}$ 
defined by $b$. If moreover $(b_1,\ldots , b_m)$ is an orthonormal basis of 
$V_1^{\perp} \cong V_0^{\ast}$ then the even part 
${\alpha ({\fr p}(V))}^b_0 \cong {\fr spin} (k,l)$ is a canonically embedded 
spinor Lie algebra, i.e.\ 
\[ {\alpha ({\fr p}(V))}^b_0 = {\fr spin}_{\sigma} := 
\{ \left( \begin{array}{ll} 
A & 0 \\
0 & \sigma (A)
\end{array} \right) | A\in {\fr so} (k,l) \subset {\fr gl}_m\} \, ,\]
where $\sigma : {\fr so} (k,l) \rightarrow {\fr gl}_n$ is equivalent
to the spinor representation. 

The linear group ${Spin}_{\sigma}\subset GL_{m|n}(\mbox{\R})$ generated 
by the Lie algebra ${\fr spin}_{\sigma} \subset ({\fr gl}_{m|n})_0$  
$\cong {\fr gl}_m \oplus{\fr gl}_n$ acts on the set of bases of $V^{\ast}$ 
from the right. 

\begin{prop} \label{b1Prop} 
Assume that $\alpha ({\fr p}(V))^b_0 = {\fr spin}_{\sigma}$
and $b'  = bg$ for some $g\in {Spin}_{\sigma}$. Then $\alpha ({\fr p}(V))^b
= \alpha ({\fr p}(V))^{b'}$.
\end{prop} 

\noindent
{\bf Proof:} This follows from the fact that 
$\alpha ({\fr p}(V))^b_0 = {\fr spin}_{\sigma}$
and $\alpha ({\fr p}(V))^b_1 = \alpha (V_1)^b$ are invariant under
${\fr spin}_{\sigma} = \alpha ({\fr spin} (V_0))^b$. $\Box$ 

Now let $(e_1,\ldots ,e_m)$ be an orthonormal basis of $V_0$ and 
$(\theta^1,\ldots ,\theta^n)$ a basis of $V_1$. The dual bases of $V_0^{\ast}$
and $V_1^{\ast}$ will be denoted by $(e^i)$ and $(\theta_j)$. 

\begin{prop} \label{b2Prop} With respect to the basis 
$b = (e^1,\ldots ,e^m,\theta_1,\ldots ,\theta_n)$ of $V^{\ast} 
\cong V_0^{\ast} + V_1^{\ast}$ the super Lie algebra 
$\alpha ({\fr p}(V)) \subset {\fr gl} (V^{\ast})$ is identified with
\[ \alpha ({\fr p}(V))^b = 
\{ \left( \begin{array}{ll} 
A & 0 \\
C & \sigma (A)
\end{array} \right)| A\in {\fr so} (k,l) ,\: C^{ji} = e^i (\pi (s \vee
\theta^j)), \: s\in V_1\}   \, ,\]           
where $C = (C^{ji})$, $j = 1,\ldots , n$, $i = 1,\ldots , m$, and $\sigma :
{\fr so} (k,l) \rightarrow {\fr gl}_n$ is equivalent
to the spinor representation. 
\end{prop}

\bigskip\noindent
2.3 {\bf The (pseudo) Riemannian supergeometry associated to the spinor 
bundle.} Now we carry over the construction of 2.2 to the 
$\mbox{\Z}_2$-graded vector bundle $V := TM_0 + S$ over $M_0$. We assume that
$M_0$ is simply connected. 
The vector bundle $V$ carries the canonical connection induced by the 
Levi-Civita connection of the pseudo Riemannian manifold $(M_0,g_0)$. 
The holonomy algebra of $V$ at $p\in M_0$ is a subalgebra of ${\fr spin} 
(T_pM_0)\subset {\fr gl} (V_p)_0$. This implies, in particular, that the 
bundle of ${Spin} (TM_0)$-invariant bilinear forms on $S$ is flat. Let 
$g_1$ be a parallel non-degenerate suitable bilinear form on $S$, see Def.\ 
\ref{suitableDef} and the remarks following Thm.\ 1. 

The ${Spin} (TM_0)$-invariant bilinear form $g = g_0 + g_1$ on $V$ should be 
thought of as a pseudo Riemannian metric for the supermanifold $M = M(S)$. 
Note that, due to Prop.\ \ref{T_pMProp}, $g(p)$ induces 
a non-degenerate bilinear form on $T_pM$. 
However, recall that $g_1$ is 
symmetric or skew-symmetric. The map $\pi = \pi _{\rho ,g_1} : 
\vee^2 S \rightarrow TM_0$ defines on ${\fr p} (V) = {\fr spin} (TM_0) + 
S\subset {\fr gl} (V)$ the structure of bundle of Poincar\'e super Lie 
algebras. ${\fr p} (V)$ is a parallel bundle. Now let 
$\alpha : {\fr p} (V) \rightarrow {\fr gl} (V^{\ast})$ be the field of 
representations induced by the coadjoint representation, cf.\ Prop.\ 
\ref{coadProp}. The image $\alpha ({\fr p} (V)) \subset {\fr gl} (V^{\ast})$
is a parallel bundle of super Lie algebras. 

\begin{prop} The frame bundle of $V^{\ast} \rightarrow M$ has a subbundle
$P_{{Spin}_{\sigma}}$ with structure group ${Spin}_{\sigma} \subset 
GL_{m|n} (\mbox{\R})$, ${Spin}_{\sigma} \cong {Spin} (k,l)$, such that 
for all $b = (e^i,\theta_j) \in (P_{{Spin}_{\sigma}})_p$: 
\begin{enumerate}
\item[1)] $(e^i)$ is an orthonormal basis of $T_p^{\ast} M_0$ and 
\item[2 )] $\alpha ({\fr p} (V_p))$ is identified via $b$ with the subalgebra
${\fr g} = \alpha ({\fr p} (V_p))^b \subset {\fr gl}_{m|n} (\mbox{\R})$,
where
\[ {\fr g}_0 = {\fr spin}_{\sigma} = \{ \left( \begin{array}{ll} 
A & 0 \\
0 & \sigma (A)
\end{array} \right)| A\in {\fr so} (k,l) \} \quad \mbox{and} \]
\[ {\fr g}_1 = \{ \left( \begin{array}{ll} 
0 & 0 \\
C & 0
\end{array} \right)| C = (C^{ji}),\: C^{ji} = e^i (\pi (s \vee
\theta^j)), \: s\in S_p   \} \] 
are independent of $b$ and $p$. Here $(\theta^j)$ is the basis of $S_p$ dual
to  $(\theta_j)$. 
\end{enumerate} 
\end{prop}

\noindent
{\bf Proof:} This follows from the holonomy reduction and Propositions 
\ref{b1Prop} and \ref{b2Prop}. $\Box$ 

We denote by $\cal V$ the sheaf of local sections of $V$. Identifying $TM_0$ 
and $T^{\ast}M_0$ via $g_0$, the map $\iota$ of Prop.\ \ref{iotaProp} 
corresponds to a monomorphism 
$\iota : {\cal V} = {\cal T}^{\ast}_{M_0} + {\cal S}^{\ast} \hookrightarrow
{\cal T}_M$. This induces a map 
\[ \iota : \Gamma (U, P_{{Spin}_{\sigma}}) 
\rightarrow {\cal F} (U)\, ,\]
where ${\cal F} (U)$ is the set of frame fields of $M$ over the open set 
$U\subset M_0$.  The image of $\iota$ generates a ${Spin}_{\sigma}$-structure
on $M$, where ${Spin}_{\sigma}$ is now considered as (purely even) linear 
supergroup ${Spin}_{\sigma} \subset GL_{m|n}$. More precisely, 
recall that ${Spin}_{\sigma} ({\cal A}(U))$ is the group generated by 
$\exp {\fr spin}_{\sigma} ({\cal A}(U)) \subset GL_{m|n} ({\cal A}(U))$. It 
acts on ${\cal F} (U)$ from the right. 
Put 
\[ {\cal F}_{{Spin}_{\sigma}} (U) \: := \: 
\iota (\Gamma (U,P_{{Spin}_{\sigma}})) {Spin}_{\sigma}({\cal A}(U))\, .\]

\begin{prop} ${\cal F}_{{Spin}_{\sigma}}$ is a ${Spin}_{\sigma}$-structure
on $M$. 
\end{prop}

Denote by $G$ the linear supergroup defined by the linear super Lie algebra
$\fr g$, see Example 3. Since ${\fr spin}_{\sigma} \subset {\fr g}
\subset {\fr gl}_{m|n} (\mbox{\R})$, we have the following inclusions of 
linear supergroups:
\begin{equation} {Spin}_{\sigma} \subset G \subset GL_{m|n} \, .
\end{equation} 
Put ${\cal F}_G (U) := {\cal F}_{{Spin}_{\sigma}} (U) G({\cal A}(U))$ for all
open $U \subset M_0$. 

\begin{prop} ${\cal F}_G$ is a $G$-structure on $M$.
\end{prop}

\begin{dof} \label{KillingvectorDef} A {\bf Killing vector field} on 
$(M,g)$ is an infinitesimal automorphism of the $G$-structure ${\cal F}_G$,
see Def.\ \ref{infautoDef}. 
\end{dof} 

\bigskip\noindent
2.4 {\bf Twistor spinors as Killing vector fields.} 
\begin{dof} \label{twistorDef} A section $s$ of the spinor bundle 
$S\rightarrow M_0$ is called
a {\bf twistor spinor} if there exists a section $\tilde{s}$ of $S$ such that
\begin{equation} \label{twistorEqu} \nabla_X s = \rho (X) \tilde{s} 
\end{equation} 
for all vector fields $X$ on $M_0$. Here $\rho (X) : S \rightarrow S$ is 
Clifford multiplication. A twistor spinor $s$ is called a Killing spinor
if $\tilde{s} = \lambda s$ for some constant $\lambda \in \mbox{\R}$
\end{dof}

\noindent
{\bf Remark:}  From (\ref{twistorEqu}) it follows that $\tilde{s} 
= -\frac{1}{m} Ds$, where $D$ is the Dirac operator.

\medskip\noindent 
The non-degenerate bilinear form $g_1$ on $S$ induces the isomorphism 
\[ S \stackrel{\sim}{\rightarrow} S^{\ast}, \quad s\mapsto s^{\ast}:=
g_1 (s,\cdot )\, .\] 
Recall that $\iota | {\cal S}^{\ast} : {\cal S}^{\ast}\hookrightarrow 
{\cal T}_M$ is simply given by interior multiplication, s.\ 2.1.
To any spinor field $S$ we associate the odd vector field
 $X_s := \iota (s^{\ast})$ on $M$. 
Now we can state the main result of this paper.

\begin{thm} Let $(M_0,g_0)$ be a pseudo Riemannian spin manifold with 
spinor bundle $(S,g_1)$; $g_1$ a parallel non-degenerate suitable
bilinear form on $S$, see Def.\ \ref{suitableDef} and 2.3.  Consider the 
supermanifold $M = M(S)$ with
the bilinear form $g = g_0 + g_1$ and let 
$s$ be a section of $S$. The vector field 
$X_s$ is a Killing vector field on $(M,g)$ iff $s$ is a twistor spinor, 
see Def.\ \ref{KillingvectorDef} and \ref{twistorDef}.  
\end{thm}  

\begin{cor} A Killing vector field $X_s$ for an extension $g$ of $g_0$
is a Killing vector field for any other extension; the extensions beeing
as in 2.3.  
\end{cor}      

\begin{lemma} \label{1Lemma} 
For all sections $s^{\ast}$, $t^{\ast}$ of $S^{\ast}$ and 
$X$ of $TM_0$ we have:
\begin{enumerate}
\item[i)] $[\iota (s^{\ast}) , \iota (t^{\ast})] = 0$,
\item[ii)] $[\iota (s^{\ast}) , \iota (X)] = [\iota (s^{\ast}),\nabla_X] =
-\iota ((\nabla_X)^{\ast})$.
\end{enumerate} 
\end{lemma} 

\noindent
{\bf Proof:} i) By definition of the supercommutator $[\cdot , \cdot ]$ on
${\cal T}_M$, we have $[\iota (s^{\ast}) , \iota (t^{\ast})] = 
\iota (s^{\ast})\circ \iota (t^{\ast}) + 
\iota (t^{\ast})\circ \iota (s^{\ast}) = 0$. \\
ii) Recall that $s^{\ast} = g_1 (s, \cdot )$. If $t$ is a section of $S$ we 
have $[\iota (s^{\ast}) , \iota (X)](t) = s^{\ast}(\nabla_Xt) -
\nabla_X s^{\ast}(t) = g_1 (s, \nabla_Xt) - \nabla_X g_1(s,t) = 
-g_1 (\nabla_Xs,t) = -(\nabla_Xs)^{\ast}(t)$. $\Box$ 

\begin{prop}  Let $s$ be a twistor spinor. For all vector 
fields $X$ and 
spinor fields $t$ on $M_0$ we have:
\begin{enumerate} 
\item[i)] $[\iota (s^{\ast}) , \iota (X)] = 
-\iota ( (\rho (X) \tilde{s})^{\ast}) 
= - \tau (g_1) \iota ( {\rho (X)}^{\ast} {\tilde{s}}^{\ast})$, where
$\tau (g_1) \in \{ \pm 1\}$ is the type of $g_1$, see Def.\ \ref{suitableDef}.
\item[ii)] $[\iota (s^{\ast}) , \iota (X)](t) = - g_1 (\rho (X) \tilde{s},t) 
= - g_0 (\pi (\tilde{s}\vee t),X)$. 
\end{enumerate} 
\end{prop} 

\noindent
{\bf Proof:} The first equation of i) follows from Lemma 
\ref{1Lemma} ii),
since $\nabla_X s = \rho (X) \tilde{s}$. Now the second equation of i)
and the first equation of ii) follow from the definition of the type 
$\tau$: $(\rho (X) \tilde{s})^{\ast} (t) = g_1 (\rho (X) 
\tilde{s},t) =
\tau (g_1) g_1(\tilde{s}, \rho (X) t)$. The last equation of ii)
is simply the
definition of $\pi = \pi_{\rho ,g_1}$, cf.\ (\ref{piEqu}). $\Box$

\medskip\noindent
{\bf Proof (of Theorem 2):} Let $(e^i,\theta_j) \in 
\Gamma (U,P_{{Spin}_{\sigma}})$, $U \subset M_0$ open, and 
$(e_i,\theta^j)$ the
dual local frame for $V = TM_0 + S$. Put 
\[ E := (\iota (e^i), \iota (\theta_j)) \in 
\Gamma (U,{\cal F}_{{Spin}_{\sigma}}) \subset \Gamma (U,{\cal F}_G)\, .\]
Since $(e_i)$ is orthonormal, i.e.\ $g_0(e_i,e_j) = \varepsilon_i
\delta_{ij}$, $\varepsilon_i \in \{\pm 1 \}$, we have 
$e^i = \varepsilon_i g_0(e_i,\cdot )$. Hence, by definition
of $\iota$ on ${\cal T}_{M_0}^{\ast}$, we have $\iota (e^i) =
\varepsilon_i \iota (e_i)$. Therefore by Lemma \ref{1Lemma} 
for any $s\in \Gamma (U,S)$ we have
\begin{equation} \label{1Equ} L_{X_s} E = ([X_s,\iota (e^i)], 
[X_s,\iota (\theta_j)]) =
(-\varepsilon_i \iota ((\nabla_{e_i} s)^{\ast} ), 0)\, ,
\end{equation}
\begin{equation} \label{2Equ} (\nabla_{e_i} s)^{\ast}  (\theta^j) = 
g_1 (\nabla_{e_i}s,\theta^j)\, .
\end{equation} 
{}From this computation it follows that $L_{X_s} E \in E ({\fr g} \otimes 
{\cal A} (U))$ iff there exists a $t\in \Gamma (U,S)$ such that
\begin{equation} \label{3Equ} L_{X_s} E = E C_t\, ,\quad \mbox{where} 
\end{equation}  
\begin{equation} \label{4Equ} C_t = \left( \begin{array}{ll} 0 & 0\\
                                (C^{ji}_t) & 0
\end{array} \right) \in  {\fr g} \otimes {\cal A} (U)\, , \quad 
C^{ji}_t = e^i (\pi (t \vee \theta^j))\, ,
\end{equation} 
see Prop.\ \ref{b2Prop}. By (\ref{1Equ}), (\ref{2Equ}) and (\ref{4Equ}) 
equation (\ref{3Equ}) is equivalent to
\begin{equation} \label{5Equ} g_1 (\nabla_{e_i}s,\theta^j) =
- \varepsilon_i  
e^i (\pi (t \vee \theta^j))\, , \quad i = 1,\ldots , m , \: j = 1,\ldots , n.
\end{equation}     
The right-hand-side is 
\begin{equation} \label{6Equ} - \varepsilon_i 
e^i (\pi (t \vee \theta^j)) = 
- g_0 (\pi (t \vee \theta^j), e_i) = 
- g_1 (\rho (e_i)t,\theta^j)\, ,
\end{equation}  
hence (\ref{5Equ}) is equivalent to the twistor equation (\ref{twistorEqu}) 
with $\tilde{s} = -t$. $\Box$

\medskip\noindent
{\bf Acknowledgments.} D.V.A.\ is grateful to  Max-Planck-Institut f\"ur 
Mathematik and V.C.\ thanks  
the Mathematical Sciences Research Institute, 
S.-S.\ Chern and R.\ Osserman, for hospitality and support.
V.C.\ would also like to thank W.\ Kramer, E.\ Poletaeva and V.\ Serganova 
for discussions related to the subject of this paper. Three of us 
(D.V.A., C..D., U.S.) should like to
thank the Mathematisches Forschungsinstitut Oberwolfach,
where this work was begun, for hospitality in the
framework of the Research-in-Pairs programme supported by the 
Volkswagen Stiftung.

\goodbreak

\begin{thebibliography}{A-C}
\bibitem[A-C]{A-C}  D.V.\ Alekseevsky, V.\ Cort\'es: {\it Classification
of $N$-(super)-extended Poincar\'e algebras and bilinear invariants
of the spinor representation of $Spin(p,q)$}, Commun.\ Math.\ Phys., to appear.
\bibitem[Ba]{Ba} M.\ Batchelor: {\it The structure of supermanifolds},
Trans.\ Amer.\ Math.\ Soc.\ {\bf 253} (1979), 329-338. 
\bibitem[Be]{Be} F.A.\ Berezin: {\it Introduction to superanalysis}, MPAM, 
D.\ Reidel, Dordrecht, 1987.
\bibitem[L]{L} D.A.\ Leites: {\it Introduction to the theory of 
supermanifolds}, Russian Math.\ Surveys {\bf 35} (1980), 3--57. 
\bibitem[M]{M} Yu.I.\ Manin: {\it Gauge field theory and complex geometry},
Springer-Verlag, Berlin-Heidelberg, 1988.
\bibitem[K]{K} B.\ Kostant: {\it Graded manifolds, graded Lie theory and
prequantization} in Springer LNM {\bf 570} (1977), 177--306.
\bibitem[S-W]{S-W} S.\ Shnider, R.O.\ Wells: {\it Supermanifolds, 
super twistor spaces and super Yang-Mills fields}, Presses de l'Universit\'e 
de Montr\'eal, 1989.    
\end{thebibliography}
\end{document}